# An Efficient Calculation of Quaternion Correlation of Signals and Color Images

Artyom M. Grigoryan[a], Sos S. Agaian[b]
[a]The University of Texas at San Antonio, [b]City University of New York / CSI
amgrigoryan@utsa.edu, sos.agaian@csi.cuny.edu

**ABSTRACT**

Over the past century, a correlation has been an essential mathematical technique utilized in engineering sciences, including practically every signal/image processing field. This paper describes an effective method of calculating the correlation function of signals and color images in quaternion algebra. We propose using the quaternions with a commutative multiplication operation and defining the corresponding correlation function in this arithmetic. The correlation between quaternion signals and images can be calculated by multiplying two quaternion DFTs of signals and images. The complexity of the correlation of color images is three times higher than in complex algebra.

## 1. INTRODUCTION

Correlation is one of the most fundamental concepts used in almost any engineering science. Correlation measures the degree to which two variables are related and expresses in quantitative terms the strength and direction of the relationship between these variables. Practically, correlation, autocorrelation, and phase-correlation become basic tools in many applications, including computer vision, robotics, sport, biomedicine, acoustic, human movement and rehabilitation research, quantum computation, geophysical applications, and image and signal processing applications, to name only a few [1]-[4]. For example, cross-correlation is a key statistical method for getting the degree of relationship/similarity among two templates/images [1] or analyses in the human movement and rehabilitation sciences [2].

Quaternions, an extension of complex numbers, were launched in the nineteenth century by Gauss. However, quaternion information processing has only recently become popular in engineering. Recently, quaternion algebra has become a very effective tool for color image application and for describing rotations and orientations in 3D spaces [5]-[9]. For instance, (i) many engineering challenges present themselves as three-dimensional vectors, ideal for quaternion representation, and (ii) quaternions naturally use inter-channel correlations in these three-dimensional problems. For example, in the RGB, CMY, and XYZ color models, color images with grayscale components can be represented as quaternion images. Unlike commonly used image processing applications, where color images deal with all color channels separately, quaternion algebra allows them to be one unit [10,12]. Also, many signal processing procedures have been expanded to the quaternion domain, including quaternion Fourier transforms, wavelets, Kalman filtering, color image restoration and enhancement, regression, quaternion neural networks, quaternion moments, and least mean square adaptive filtering [10]-[19]. The quaternion image processing tools are very effective in some cases, especially in the frequency domain. See, for example, the 2-D quaternion discrete Fourier transforms-based alpha-rooting method for image enhancement [11]-[13].

It is natural to extend the correlation procedure to the quaternion domain. The challenge here is that traditional quaternion algebra, which is noncommutative, does not provide effective methods of calculating convolution and correlation. Because of the noncommutativity of the multiplication, the correlation of images has two forms; when multiplying images in the double sum from left and right; the results are different. Also, the concept of quaternion discrete Fourier transforms (QDFT) is not unique, and the question arises of how to choose the transforms [25,26].

To overcome these limitations, we consider the 4-D model of signals and color images in the commutative quaternion algebra, called the (2,2)-model of quaternions, and describe an effective method of calculating the correlation function for signals and images in this paper. The key paper contributions are

1. We propose an effective tool to calculate the correlation of signals and color images in the quaternion algebra domain.
2. In the frequency-domain, the quaternion correlation of two images can be calculating by multiplying their QDFTs.
3. We estimate the complexity of the correlation of color images and show that
   a. The complexity of calculating the correlation function based on quaternion algebra is two times higher than in the case of complex algebra.
   b. Calculating the correlation function in the frequency domain using the commutative (2,2)-model can be accomplished without quaternion DFTs.
   c. The quaternion correlation needs only six complex 2-D DFTs plus $4NM$ complex multiplications or $16NM$ real multiplications.

The rest of the paper is organized in the following way. In Section 2, a brief introduction to quaternion numbers is given in the (1,3)- and (2,2)-models. In Section 3, the quaternion correlation of signals is described with an example. The quaternion correlation implementation block-diagram is illustrated. Section 4 presents the correlation of color images in the (2,2)-model. The complexity of the quaternion correlation is compared with the known methods in the traditional commutative (1,3)-model of quaternions. Finally, Section 5 presents the conclusion and future work.

## 2. QUATERNION NUMBERS

This section briefly introduces quaternion numbers in the noncommutative quaternion algebra and the (2,2)-model [23]. This material and subsequent sections will allow us to work with the concept of the correlation function of quaternion images. To simplify the discussion of this topic, we will first work on the one-dimensional case and then generalize the results to the two-dimensional case. As we know, fast correlation methods are based on the fast Fourier transform since the correlation can be represented as a linear convolution and then reduced the latter to a cyclic convolution. It should be recalled that the linear convolution and DFT do not have unique definitions in traditional quaternion arithmetic because it is not commutative. All these obstacles could be removed if we had commutative multiplication. Therefore, we first describe these concepts in commutative quantum arithmetic, which we call the (2,2)-mode, and then introduce the solution to the problem of calculating the quaternion correlation.

### 2.1. The (1,3)-model of quaternions

The traditional quaternion arithmetic was first described by Gauss [20], and later by Hamilton [21,22]. The quaternion is the two-component number $q = a + (bi + cj + dk)$. This number presents a quarter of real numbers $(a, b, c, d)$. The real part of the quaternion is $a$ and the imaginary part is the three-component number $q' = (bi + cj + dk)$. This is why, we call such a representation of quaternions the (1,3)-model [23]. The imaginary units $i, j,$ and $k$ are orthogonal to each other and to real unit 1. The multiplication laws for these imaginary units are the following:

$$ij = -ji = k, \quad jk = -kj = i, \quad ki = -ik = j, \quad i^2 = j^2 = k^2 = ijk = -1. \tag{1}$$

Thus, the quaternion number presents the 4-D vector $(a, q')$, where 3-D vector is $q' = (b, c, d)$. The conjugate of $q$ is defined as $\bar{q} = a - (bi + cj + dk)$. The length, or modulus of the quaternion is $|q| = \sqrt{q\bar{q}} = \sqrt{a^2 + b^2 + c^2 + d^2}$. In polar form, the quaternion is written as $q = |q|\exp(\mu\vartheta)$. Here, $\vartheta$ is an angle and $\mu$ is a unit pure quaternion, such that $\mu^2 = -1$ and $|\mu| = 1$. Therefore, the quaternion exponent can be calculated as

$$\exp(\mu\vartheta) = \cos(\vartheta) + \mu\sin(\vartheta).$$

Because of inequality $ij \neq ji$ and $jk \neq kj$, the multiplication is not a commutative operation in the (1,3)-model. It also means, that the fundamental multiplicative identity does not hold for quaternion exponents, i.e.,

$$\exp(\mu\vartheta_1)\exp(\mu\vartheta_2) \neq \exp(\mu[\vartheta_1 + \vartheta_2]), \text{ for all } \vartheta_1, \vartheta_2. \tag{2}$$

For some cases, this property holds. For instance, when $\mu = i, j$, and $k$.

The multiplication of two quaternions $q_1 = a_1 + q'_1 = a_1 + (b_1 i + c_1 j + d_1 k)$ and $q_2 = a_2 + q'_2 = a_2 + (b_2 i + c_2 j + d_2 k)$ is calculated by

$$q = q_1 q_2 = [a_1 q'_2 + a_2 q'_1] + a_1 a_2 - [b_1 b_2 + c_1 c_2 + d_1 d_2] + \begin{vmatrix} i & j & k \\ b_1 & c_1 & d_1 \\ b_2 & c_2 & d_2 \end{vmatrix}. \tag{3}$$

In matrix form, this quaternion $q = a + (bi + cj + dk)$ can be calculated by

$$q = \begin{bmatrix} a \\ b \\ c \\ d \end{bmatrix} = \mathbf{A}_{ij} \begin{bmatrix} a_2 \\ b_2 \\ c_2 \\ d_2 \end{bmatrix} = \begin{bmatrix} a_1 & -b_1 & -c_1 & -d_1 \\ b_1 & a_1 & -d_1 & c_1 \\ c_1 & d_1 & a_1 & -b_1 \\ d_1 & -c_1 & b_1 & a_1 \end{bmatrix} \begin{bmatrix} a_2 \\ b_2 \\ c_2 \\ d_2 \end{bmatrix}, \quad \det \mathbf{A}_{ij} = |q_1|^4. \tag{4}$$

Also, we can define the following multiplication laws:

$$ji = -ij = k, \quad jk = -kj = -i, \quad ki = -ik = -j, \quad i^2 = j^2 = k^2 = jik = -1. \tag{5}$$

With such rules, the product of two quaternions $q = q_1 q_2$ is equal to the product $q_2 q_1$ in the first arithmetic of quaternions with rules of Eq. 1. The corresponding matrix of this multiplication is equal to

$$\mathbf{A}_{ji} = \begin{bmatrix} a_1 & -b_1 & -c_1 & -d_1 \\ b_1 & a_1 & d_1 & -c_1 \\ c_1 & -d_1 & a_1 & b_1 \\ d_1 & c_1 & -b_1 & a_1 \end{bmatrix}, \quad \det \mathbf{A}_{ij} = |q_1|^4.$$

It should be noted, that these two matrices are orthogonal.

## 2.2. The (1,3)-model of quaternions

In this section, we consider the commutative (2,2)-model of quaternions and compare briefly it with the noncommutative (1,3)-model, for which the multiplication rules are defined as in Eq. 5. The quaternion is a pair of two complex numbers

$$q = a_1 + ja_2 = (a + bi) + j(c + di) = a_1 + ja_2 = [a_1, a_2]. \tag{6}$$

Here, complex numbers $a_1 = a + bi$ and $a_2 = c + di$. The quaternion number can be written as $q = [(a, b), (c, d)]$ (assuming $ji = k$). This number represents the vector $q = (a, b, c, d)$ in the 4-D space. The following notations are used in the (2,2)-model [23]. The quaternion $q$ as a pair of two complex numbers is written as:

$$q = [a_1, a_2], \quad a_1 = (a_{1,1}, a_{1,2}), \quad a_2 = (a_{2,1}, a_{2,2}). \tag{7}$$

Here, $a_{1,1}, a_{1,2}, a_{2,1}$, and $a_{2,2}$ are real numbers. The square brackets are used for quaternions. The round brackets are used for 2-D vectors, or complex numbers. In the $a_2 = 0$ case, the quaternion $q = [a_1, 0]$ is a complex number $a_1$, and if $a_1 = (a_{1,1}, 0)$, we call $q = [a_1, 0]$ real.

Consider two quaternions $q_1 = [a_1, a_2]$ and $q_2 = [b_1, b_2]$. Here, $a_1 = (a_{1,1}, a_{1,2}), a_2 = (a_{2,1}, a_{2,2})$ and $b_1 = (b_{1,1}, b_{1,2}), b_2 = (b_{2,1}, b_{2,2})$. The following operations are defined on quaternions.

a. The sum of two quaternions is defined as the component-wise operation. Thus,

$$q_1 + q_2 = [a_1, a_2] + [b_1, b_2] = [a_1 + b_1, a_2 + b_2].$$

b. The multiplication of two quaternions is calculated by

$$q_1 q_2 = [a_1, a_2][b_1, b_2] \triangleq [a_1 b_1 - a_2 b_2, a_1 b_2 + a_2 b_1]. \tag{8}$$

Consider the following four quaternion units:

1. $e_1 = [(1,0), (0,0)]$, or shortly $e_1 = (1,0) = 1$.
2. $e_2 = [(0,1), (0,0)]$,
3. $e_3 = [(0,0), (1,0)]$
4. $e_4 = [(0,0), (0,1)]$.

The multiplications of these quaternions are given in Table 1, which we call $T(1, e_2, e_3, e_4)$ table.

|     | $e_2$ | $e_3$ | $e_4$ |
| --- | --- | --- | --- |
| $e_2$ | $-1$ | $e_4$ | $-e_3$ |
| $e_3$ | $e_4$ | $-1$ | $-e_2$ |
| $e_4$ | $-e_3$ | $-e_2$ | $1$ |

Table 1. $T(1, e_2, e_3, e_4)$.

In matrix form, the product of two quaternions $q = q_1 q_2 = [a_1, a_2][b_1, b_2]$ can be written as

$$q = \begin{bmatrix} a \\ b \\ c \\ d \end{bmatrix} = M_{q_1} q_2 = \begin{bmatrix} a_{1,1} & -a_{1,2} & -a_{2,1} & a_{2,2} \\ a_{1,2} & a_{1,1} & -a_{2,2} & -a_{2,1} \\ a_{2,1} & -a_{2,2} & a_{1,1} & -a_{1,2} \\ a_{2,2} & a_{2,1} & a_{1,2} & a_{1,1} \end{bmatrix} \begin{bmatrix} b_{1,1} \\ b_{1,2} \\ b_{2,1} \\ b_{2,2} \end{bmatrix}. \tag{9}$$

One can notice the quaternion number $q_1$ in the first column of this matrix. The matrix $M_{q_1}$ is not orthogonal.

***Example 1.*** Consider two quaternion numbers $q_1 = (a_1, a_2) = [(1,4), (-1,2)]$ and $q_2 = (b_1, b_2) = [(2,5), (3,-1)]$. The product $q = q_1 q_2$ is calculated by the matrix $M$ of multiplication

$$M = M_{q_1} = \begin{bmatrix} 1 & -4 & 1 & 2 \\ 4 & 1 & -2 & 1 \\ -1 & -2 & 1 & -4 \\ 2 & -1 & 4 & 1 \end{bmatrix}, \quad \det M = 340.$$

The inverse matrix $M^{-1}$ as

$$M^{-1} = \frac{1}{170}\begin{bmatrix} -1 & 38 & 13 & 16 \\ -38 & -1 & -16 & 13 \\ -13 & -16 & -1 & 38 \\ 16 & -13 & -38 & -1 \end{bmatrix}.$$

The $q = q_1 q_2$ is calculated by

$$q = \begin{bmatrix} 1 & -4 & 1 & 2 \\ 4 & 1 & -2 & 1 \\ -1 & -2 & 1 & -4 \\ 2 & -1 & 4 & 1 \end{bmatrix} \begin{bmatrix} 2 \\ 5 \\ 3 \\ -1 \end{bmatrix} = \begin{bmatrix} -17 \\ 6 \\ -5 \\ 10 \end{bmatrix}.$$

Therefore, $q = q_1 q_2 = [(-17,6), (-5,10)]$. We note for comparison that in the (1,3)-model, the product of these two quaternions is equal to $[-13,8,5,20]$ which is the quaternion $-13 + (8i + 5j + 20k)$. The results are different.

The properties of this multiplication are described in detail in work [13] of the authors. Here, we briefly mention the main properties:

1. If $k$ is a real or complex number, then $kq = (ka_1, ka_2)$.
2. The multiplication is associative $(q_1 q_2)q_3 = q_1(q_2 q_3)$, for any quaternions $q_1$, $q_2$, and $q_3$.
3. The multiplication is distributive

$$q_1(q_2 + q_3) = q_1 q_2 + q_1 q_3, \qquad (10)$$

for any quaternions $q_1$, $q_2$, and $q_3$
4. The multiplication is commutative,

$$q_1 q_2 = q_2 q_1.$$

5. The square of the quaternion is

$$q^2 = qq = [a_1^2 - a_2^2, 2a_1 a_2]. \qquad (11)$$

6. The zero has divisors, i.e., exist such quaternions $q_1$ and $q_2 \neq 0$, that $q_1 q_2 = 0$. For instance,

$$(1 + e_4)(1 - e_4) = 1 - e_4^2 = 0.$$

Thus, $(1 + e_4)$ and $(1 - e_4)$ are the divisors of the zero.
7. The conjugate of the quaternion $q$ is the number

$$\bar{q} = [\bar{a}_1, \bar{a}_2] = [(a_{1,1}, -a_{1,2}), (a_{2,1}, -a_{2,2})]. \qquad (12)$$

The conjugates of unit quaternions $\bar{e}_2 = [(0,-1),(0,0)] = -e_2$, $\bar{e}_3 = e_3$, and $\bar{e}_4 = [(0,0),(0,-1)] = -e_4$,
8. The multiplication of a quaternion on its conjugate is equal to

$$q\bar{q} = [a_1 \bar{a}_1 - a_2 \bar{a}_2, a_1 \bar{a}_2 + a_2 \bar{a}_1] = [|a_1|^2 - |a_2|^2, a_1 \bar{a}_2 + a_2 \bar{a}_1]. \qquad (13)$$

This number is not the square of the modulus of the quaternion, as in the non-commutative (1,3)-model. However, when a quaternion is a complex number, $q = [a_1, 0]$, then, $q\bar{q} = |a_1|^2$. The module of the quaternion is calculated as $|q| = \sqrt{|a_1|^2 + |a_2|^2}$.

## 3. QUATERNION CORRELATION IN THE (2,2)-MODEL

This section considers the operation of correlation in the (2,2)-model of quaternions. We first describe this operation in the 1-D case of signals. The 2-D case of images is considered similarly. Let $v_k = [v_{1,k}, v_{2,k}]$, $k = 0:(L-1)$, be the quaternion signal which will be correlated with another signal $q_n = [f_n, g_n]$, $n = 0:(N-1)$. We consider that $L \leq N$.
The correlation of these signals is defined as

$$r_n(v, q) = r_n = \sum_k v_{k-n} q_k, \quad n = -(L-1):(N-1). \tag{14}$$

This operation can be written as the convolution

$$r_n(v, q) = \hat{v}_n * q_n = \sum_k \hat{v}_{n-k} q_k, \tag{15}$$

where the quaternion signal $\hat{v}_n$ is the time-reversal transform $\hat{v}_n = v_{-n}$. The linear convolution can be reduced to the aperiodic convolution. Without loss of generality, we consider that both quaternion signals $q_n$ and $v_n$ are periodic and have the same length $N$. The convolution is unique in the (2,2)-model, and the multiplicative property holds. In other words, the aperiodic quaternion convolution can be calculated by multiplication of the $N$-point QDFTs of signals [23].

The aperiodic convolution $R_n = \hat{v}_n * q_n$ can be written as follows (see more detail in [SP]):

$$r_n = [r_{1,n}, r_{2,n}] = [\hat{v}_{1,n} * f_n - \hat{v}_{2,n} * g_n, \hat{v}_{1,n} * g_n + \hat{v}_{2,n} * f_n]. \tag{16}$$

This equation can be written in term of correlations of components of the quaternion signals,

$$r_n = [r_n(v_1, f) - r_n(v_2, g), r_n(v_1, g) + r_n(v_2, f)]. \tag{17}$$

The cross-correlation functions of components of both signals are calculated, and then the sum and difference of the mixed correlations are calculated. Thus, to calculate the quaternion correlation, four aperiodic convolutions in complex arithmetic can be used, as shown in Fig. 1

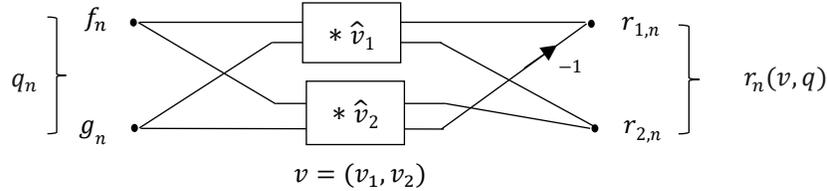

Figure 1. The diagram of the quaternion correlation in the (2,2)-model.

In the case when the quaternion signals are equal, the autocorrelation function is calculated by

$$r_n = [r_n(f, f) - r_n(g, g), r_n(f, g) + r_n(g, f)]. \tag{18}$$

When the signals $q$ and $v$ are real, i.e., $g = 0, v_2 = 0$, and $f$ and $v_1$ are real signals, the correlation $r_n(v, q) = r_n(v_1, f)$. This is the traditional case and this correlation is usually normalized as

$$r_n(v_1, f) \to R_n(v_1, f) = \frac{r_n(v_1, f)}{E[v_1] E[f]},$$

where the coefficients

$$E[f] = \sqrt{r_0(f, f)} = \sqrt{\sum_{n=0}^{N-1} f_n^2} \quad \text{and} \quad E[v_1] = \sqrt{r_0(v_1, v_1)} = \sqrt{\sum_{n=0}^{N-1} v_{1,n}^2}.$$

Here, $r_0(f,f)$ and $r_0(v_1,v_1)$ are values of the autocorrelation functions at point $n=0$. The correlation of two quaternion signals $q$ and $v$ should also be normalized. Determining the normalization for quaternion convolution is not an easy task. We think to normalize the quaternion correlation as follows:

$$r_n(v,q) \to R_n(v,q) = \left[\frac{r_n(v_1,f) - r_n(v_2,g)}{K_1}, \frac{r_n(v_1,g) + r_n(v_2,f)}{K_2}\right], \tag{19}$$

where the coefficients

$$K_1 = E[f]E[v_1] + E[g]E[v_2], \quad K_2 = E[g]E[v_1] + E[f]E[v_2].$$

We also can consider the normalization of the correlation function in a traditional way

$$r_n(v,q) \to R_n(v,q) = \frac{r_n(v,q)}{E[v]E[q]}, \tag{20}$$

where $E[q] = \sqrt{E[f]^2 + E[g]^2}$ and $E[v] = \sqrt{E[v_1]^2 + E[v_2]^2}$.

***Example 1.*** Consider the quaternion signals $q_n = [f_n, g_n]$ and $v_n = [v_{1,n}, v_{2,n}]$ of length $N = 512$ each. Each signal consists of four columns of a grayscale "jetplane" image, which is shown in Fig. 2 in part (a).

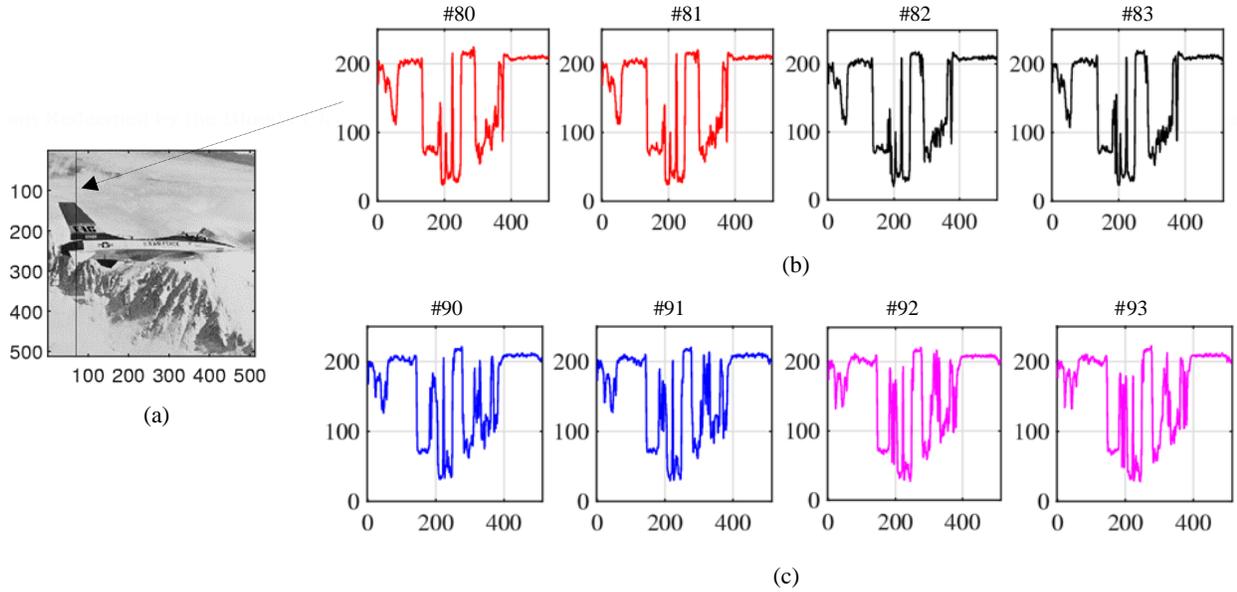

**Fig. 2** (a) The grayscale image, (b) the first quaternion signal $q$, and (c) the second quaternion signal $v$.

For signal $q_n$, four columns are taken, starting with the number 80. This signal is shown in part (b). The second signal consists of columns number 90-93 of the image, which are shown in part (c). All these columns are close to each other, so it is obvious that these two quadruples of signals are correlated. The correlation function of these two quaternion signals is shown in Fig. 3. The maximum of the quaternion correlation can be seen in the last component of the correlation at point $n = 0$. The value of this maximum is 0.9696. The amplitudes of the first three correlation components are small; they lie within the interval [-0.01,0.01].

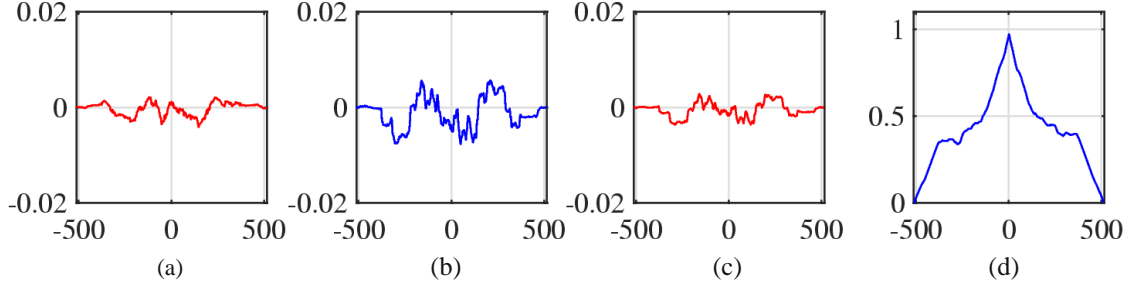

**Fig. 3** Four components of the quaternion correlation $r_n(v,q)$.

For this example, the normalized coefficients are almost the same. Indeed, $K_1 = 5.7722 \times 10^7 \approx K_2 = 57723 \times 10^7$, and the normalized quaternion correlation can be considered as

$$R_n(v,q) = r_n(v,q)\,\frac{1}{K_1}.$$

Note that the maximum correlation of two real parts of these signals, i.e., the 80th and 90th columns of the 'jetplane' image, is 0.9688 at point $n = 0$. For columns number 81 and 91 such correlation number is 0.9672, and so on. The maximum 0.9696 of quaternion correlation is close to these numbers. This fact shows that the normalization of the quaternion correlation function given in Eq. 19 can be considered successful. Perhaps there are other better normalization coefficients for signal correlation. The normalization of the correlation function in Eq. 20 results in almost the same result, since the constant $E[v]E[q] = 5.7577 \times 10^7$.

The 1021-point correlation function was calculated using the 1021-point discrete Fourier transforms of these signals after zero padding. In the frequency domain, Eq. 16 can be written as

$$R_p = [R_{1,p}, R_{2,p}] = [V_{1,N-p}F_p - V_{2,N-p}G_p, V_{1,N-p}G_p + V_{2,N-p}F_p], \qquad p = 0:(N-1). \tag{21}$$

Here, the capital letters are used for the $N$-point DFTs of the corresponding complex signals. In the (2,2)-model, this expression equals to the product $[V_{1,N-p}, V_{2,N-p}][F_p, G_p]$, which is the quaternion DFT, namely the $e_2$-QDFT (see for detail [23]). In other words, the $e_2$-QDFT of the correlation is the product of the $e_2$-QDFTs of quaternion signals $\hat{v}_n$ and $q_n$,

$$\mathcal{F}[r]_p = [V_{1,N-p}, V_{2,N-p}][F_p, G_p] = [V_{1,N-p}F_p - V_{2,N-p}G_p, V_{1,N-p}G_p + V_{2,N-p}F_p]. \tag{22}$$

The required number of operations in the frequency domain is defined by four $N$-point DFTs plus $4N$ complex multiplications and $2N$ additions. Here, the $N$-point $e_2$-QDFT of a time-reversal signal $\hat{v}_n = v_{-n}$ is calculated by

$$\mathcal{F}[\hat{v}]_p = \sum_{n=0}^{N-1} v_{-n} W_\mu^{np} = \sum_{n=0}^{N-1} v_{N-n} W_\mu^{-(N-n)p} = \sum_{n=0}^{N-1} v_n W_\mu^{-np} = \sum_{n=0}^{N-1} v_n W_\mu^{n(N-p)} = V_{N-p}, \tag{23}$$

where $p = 0:(N-1)$. The exponential function is calculated for $\mu = -e_2$ by

$$W_\mu = \exp\left(\frac{\mu 2\pi}{N}\right) = \left[\left(\cos\left(\frac{2\pi}{N}\right), -\sin\left(\frac{2\pi}{N}\right)\right), (0,0)\right] = \left[e^{-i\left(\frac{2\pi}{N}\right)}, 0\right]. \tag{24}$$

Therefore, the $N$-point QDFT of the correlation can be calculated by

$$\mathcal{F}[r]_p = \sum_{n=0}^{N-1} r_n W_\mu^{np} = V_{N-p} Q_p, \qquad p = 0\colon (N-1). \tag{25}$$

Here, $V_p$ and $Q_p$ are coefficients of the $N$-point QDFT of signals $v_n$ and $q_n$, respectively. A similar equation holds for the traditional correlation function.

### 3.1 Comparison

Now, we compare the computation of the quaternion correlation in the noncommutative algebra, or (1,3)-model, which is considered in the form

$$r_n = \sum_k q_k v_{k-n}, \qquad n = -(L-1)\colon (N-1). \tag{26}$$

Considering zero padding of both signals, the correlation $r_n$ can be calculated in the frequency domain as follows (see for detail [24]):

$$\mathcal{F}_j[r]_p = \mathcal{F}_j[q]_p \mathcal{F}_j[v_{1,1} + j v_{2,1}]_{N-p} + \mathcal{F}_j[q]_{N-p} \mathcal{F}[v_{1,2} + j v_{2,2}]_{N-p}, \tag{27}$$

where $p = 0\colon (N-1)$. Here, $N = L + N - 1$ and the $N$-point QDFT in the (1,3)-model is defined by the equation

$$\mathcal{F}_j[r]_p = \sum_{n=0}^{N-1} r_n W_\mu^{np}, \qquad p = 0\colon (N-1)$$

The exponential coefficients $W_\mu^{np}$ are defined for the quaternion basis unit $\mu = -j$. As follows from Eq. 27, the correlation requires the $N$-point QDFT of $q_n$, two DFTs of complex signals $z_1 = (v_{1,1} + j v_{2,1})$ and $z_2 = (v_{1,2} + j v_{2,2})$, and $2N$ operations of quaternion multiplication. Also, the $N$-point inverse QDFT is required. Note that the $N$-point QDFT can be calculated by two $N$-point DFTs [12]. The quaternion transforms $\mathcal{F}_j[z_1]$ and $\mathcal{F}_j[z_2]$ can be considered as the complex $N$-point DFT each. Multiplying of quaternions by complex numbers requires 8 real multiplications. Therefore, the minimum number of quaternion multiplications can be estimated as $m_{QC} = 6 m_{DFT} + 16N$. Here, $m_{DFT}$ is the number of real multiplications for the complex DFT. The diagram of the calculation of the QDFT of the correlation by Eq. 27 is shown in Fig. 4. The transform $v_n \to (z_{1,n}, z_{2,n})$ is denoted by $T$.

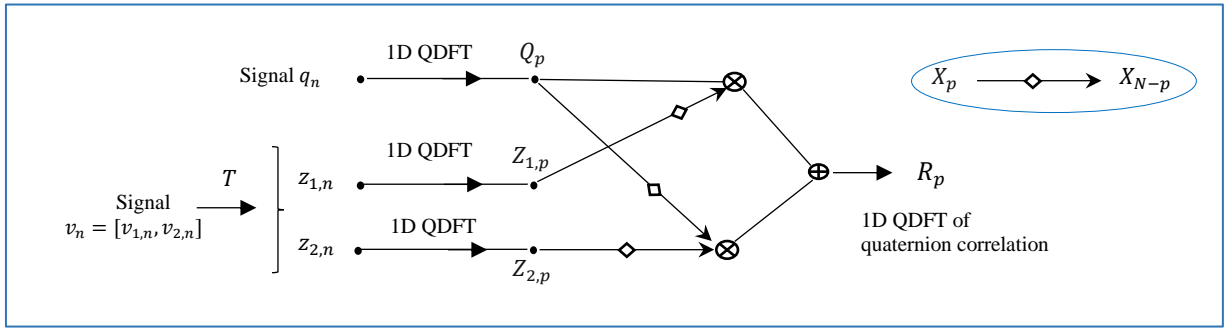

**Fig. 4** The operation of the calculation of the 1-D QDFT of the quaternion correlation $R_p$ by Eq. 27.

In the commutative (2,2)-model, the calculation of the correlation function in the frequency domain is the multiplication of the QDFTs, i.e., $R_p = V_{N-p} Q_p$. Thus, the quaternion correlation requires three QDFTs and number of quaternion multiplications is estimated as $m_{QC} = 3m_{QDFT}$. The inverse 1-D QDFT requires two complex DFTs [12]. To estimate the number of complex multiplications, we can use Eq. 21. The calculations at each frequency-point are performed by the butterfly-type operation shown in Fig. 5. Thus, the calculation of the correlation function in the frequency domain does not require quaternion DFTs, only six complex DFTs plus $4N$ complex multiplications, or $16N$ real multiplications. Thus, the estimation of complex multiplications also equals $m_{QC} = 6m_{DFT} + 16N$. Quaternion multiplication operations are not required.

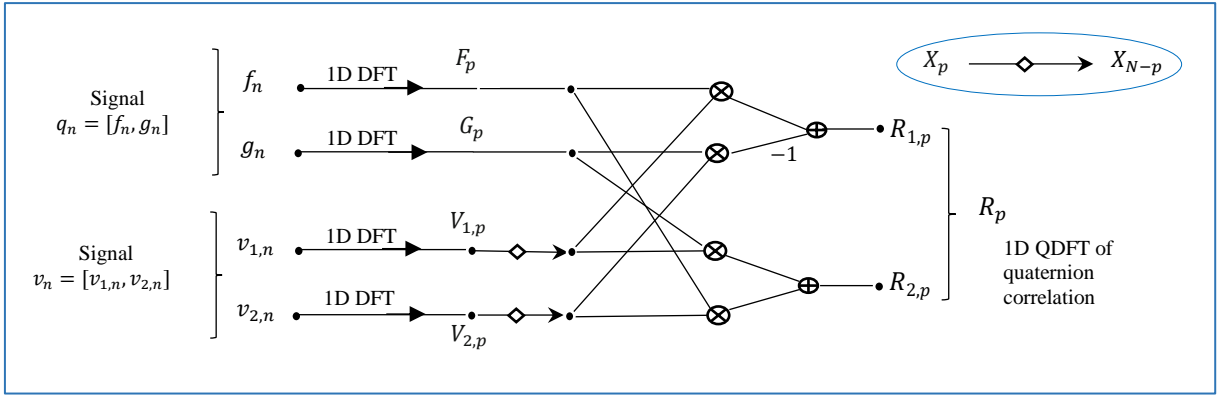

**Fig. 5** The block-diagram of calculation of the 1-D QDFT of the quaternion correlation $R_p(v, q)$.

The same operation of a butterfly is shown in Fig. 6 in a more compact form.

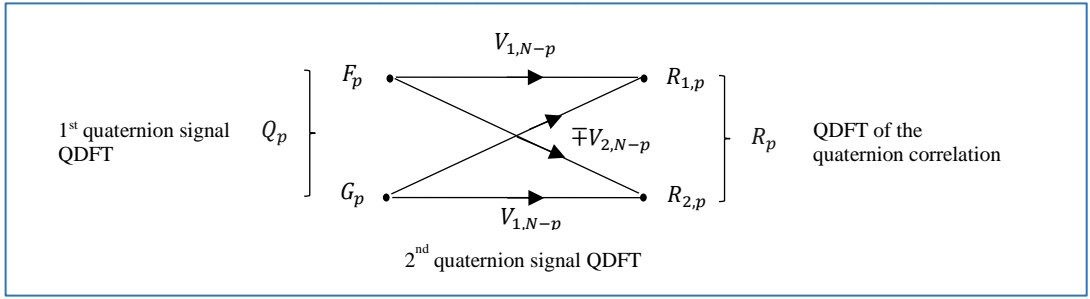

**Fig. 6** The operation of the 2×2 butterfly for the quaternion correlation $R_p(v, q)$.

The main difference between these two methods is the following:
- In the (2,2)-model, the quaternion correlation s $r_n$ of two signals $v_n$ and $q_n$ is defined by four traditional cross correlation functions of their components. No quaternion operations are required.
- Also, in the frequency domain, the correlation of two signals is the operation of multiplication of their QDFTs, $R_p = V_{N-p}Q_p, p = 0:(N-1)$.
- In the traditional (1,3)-model, this property does not hold, i.e., $R_p \neq V_{N-p}Q_p$, for either left or right QDFT.

It should be noted that the quaternion convolutions in the (1,3)- and (2,2)-models are different. The difference is interesting and significant, and it is easy to see in the above example. The correlation of the quaternion signals $q$ and $v$ in the (1,3)-model is shown in Fig. 7. This correlation was normalized as in Eq. 20. All four components show almost the same correlation (the first component is similar up to the sign) with maximum correlation in point $n = 0$. In the (2,2)-model, only the last component of the correlation is so clearly expressed.

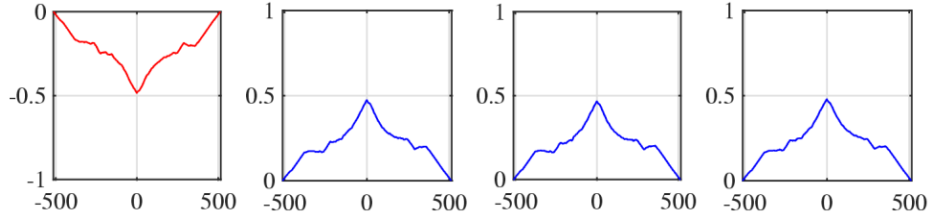

**Fig. 7** Four components of the quaternion correlation $r_n(v, q)$ in the (1,3)-model.

## 4. Quaternion 2-D Correlation of Images

Now we consider the 2-D correlation function, the correlation of two images. Let $v_{k,l} = [v_{1;k,l}, v_{2;k,l}]$ be the quaternion image of $L_1 \times L_2$ pixels. Another quaternion image $q_{n,m} = [f_{n,m}, g_{n,m}]$ is of $N_1 \times N_2$ pixels. The correlation of these images is defined as

$$r_{n,m} = r_{n,m}(v, q) = \sum_m \sum_l v_{k-n,l-m} q_{k,l}. \tag{28}$$

In the commutative (2,2)-model, this correlation function can be calculated in the frequency-domain, similar to the 2-D convolution. These two operations are reduced to the multiplication of 2-D QDFT of images. Indeed, the 2-QDFT of the correlation

$$R_{p,s} = \sum_{n=0}^{N-1} \sum_{m=0}^{M-1} r_{n,m}(v, q) W_{\mu,N}^{ms} W_{\mu,M}^{np}, \quad p = 0{:}(N-1), s = 0{:}(M-1), \tag{29}$$

is equals to $R_{p,s} = V_{N-p,M-s} Q_{p,s}$. Here, $Q_{p,s}$ and $V_{p,s}$ are the 2-D QDFT of the quaternion images $q_{n,m}$ and $v_{n,m}$, respectively. We assume zero padding the images $v_{n,m}$ and $q_{n,m}$ to the same size $N \times M$, where $N = N_1 + L_1 - 1$ and $M = N_2 + L_2 - 1$. The $N \times M$-point inverse QDFT allows us to calculate $r_{n,m}(v, q)$ from its 2-D QDFT, $R_{p,s}$. Thus, the calculation of the correlation requires three 2-D $N \times M$-point quaternion DFTs.

Also, we can consider the direct calculations of the quaternion correlation,

$$r_{n,m}(v, q) = \sum_m \sum_l v_{k-n,l-m} q_{k,l} = \sum_m \sum_l [v_{1;k-n,l-m}, v_{2;k-n,l-m}][f_{n,m}, g_{n,m}]$$
$$= \sum_m \sum_l [v_{1;k-n,l-m} f_{n,m} - v_{2;k-n,l-m} g_{n,m}, v_{1;k-n,l-m} g_{n,m} + v_{2;k-n,l-m} f_{n,m}]$$
$$= \left[\sum_m \sum_l v_{1;k-n,l-m} f_{n,m} - v_{2;k-n,l-m} g_{n,m}, \sum_m \sum_l v_{1;k-n,l-m} g_{n,m} + v_{2;k-n,l-m} f_{n,m}\right].$$

Thus, we obtain the following:

$$r_{n,m}(v, q) = [r_{n,m}(v_1, f) - r_{n,m}(v_2, g), r_{n,m}(v_1, g) + r_{n,m}(v_2, f)]. \tag{30}$$

Four traditional correlation functions are required to calculate the quaternion correlation. In the square brackets, the correlation functions are cross-correlations of components of the quaternion images. Eq. 30 can also be written in the following form with four 2-D convolutions:

$$r_{n,m}(v,q) = [(\hat{v}_1 * f)_{n,m} - (\hat{v}_2 * g)_{n,m}, (\hat{v}_1 * g)_{n,m} + (\hat{v}_2 * f)_{n,m}]. \tag{31}$$

Here, $\hat{v}_{1;n,m} = v_{1;-n,-m}$ and $\hat{v}_{2;n,m} = v_{2;-n,-m}$. In the case of the autocorrelation, i.e., when $v_{n,m} = q_{n,m} = [f_{n,m}, g_{n,m}]$, we obtain

$$r_{n,m}(q,q) = [r_{n,m}(f,f) - r_{n,m}(g,g), r_{n,m}(f,g) + r_{n,m}(g,f)],$$

or

$$r_{n,m}(q,q) = \left[(\hat{f} * f)_{n,m} - (\hat{g} * g)_{n,m}, (\hat{f} * g)_{n,m} + (\hat{g} * f)_{n,m}\right]. \tag{32}$$

Note that, $r_{n,m}(g,f) = r_{-n,-m}(f,g)$. Thus, two autocorrelations of $f$ and $g$ are required, plus the cross-correlation of these components.

The block-diagram similar to one shown in Fig. 5 can be used to calculate this 2-D quaternion correlation function. It follows from Eq. 30, that the correlation $r_{n,m}$ can be calculated by four complex $N \times M$-point 2-D DFTs of components of the quaternion images, $f_{n,m}, g_{n,m}, v_{1;n,m}$, and $v_{2;n,m}$. Also, two inverse complex $N \times M$-point 2-D DFTs. Therefore, the total number of complex operations of real multiplication for quaternion convolution in the (2,2)-model can be estimated as $m_{2DQC} = 6m_{2DDFT} + 4(4MN)$. It is assumed that the complex multiplication is performed with 4 real multiplications. Note that the traditional complex 2-D correlation requires three complex $N \times M$-point DFTs.

For comparison, we mention the calculation method of the 2-D quaternion correlation in the traditional (1,3)-model, which is described in [25]. The authors suggest calculating this operation by using the type-3 QDFT. Type-3 QDFT is the right-side transform

$$Q_{p,s} = \sum_{m=0}^{M-1}\sum_{n=0}^{N-1} q_{n,m}(v,q) e^{-\mu 2\pi[sm/M+np/N]}, \quad p = 0:(N-1), s = 0:(M-1).$$

Here, $\mu$ is a pure quaternion unit. This transform can be calculated by two complex 2-D DFT, when using the method of symplectic decomposition [26] of the imaginary part $q'_{n,m}$ of the quaternion image $q_{n,m}$. For that, the image is presented in the new basis $\{1, \mu_1, \mu_2, \mu_3\}$, where quaternion units $\mu_1$ and $\mu_2$ are orthogonal to each other, and $\mu_3 = \mu_1\mu_2$ (see detail in [12]). In general, except the cases like $\mu_1 = i$, such decomposition requires 9 multiplications for each pixel. Or, total $9NM$ real multiplications. The quaternion correlation was described for the 2-D functions, not discrete images. If we transform Eq. 115 in [26] into the discrete case, we obtain the following equation:

$$R_{p,s} = Q_{p,s}\overline{(V_1)}_{p,s} - Q_{N-p,M-s} \cdot j(V_2)_{N-p,M-s}. \tag{33}$$

Here, $(V_1)_{p,s}$ is the 2-D QDFT of the complex data $v_{1;n,m}$ and $(V_b)_{p,s}$ is the 2-D QDFT of the complex data $(v_b)_{n,m} = v_{1;n,m}$. Thus, for each frequency-point $(p, s)$, two operations of quaternion multiplications are required. The total number of multiplications, can be estimated as $m_{2DQC} \geq 6m_{2DDFT} + k2NM$. Here, $k$ stands for the number of real multiplications for one quaternion multiplication. We consider $k = 8$, as in the 1-D case above.

Summarize the result of this comparison, we can state the following (see also Table 1):
- In the (2,2)-model,
    1. The 2-D quaternion correlation s $r_{n,m}$ of two images $v_{n,m}$ and $q_{n,m}$ is defined by four traditional cross correlation functions of their components. No quaternion operations are required.

2. In the frequency-domain, the 2-D quaternion correlation of two images is the operation of multiplication of their 2-D QDFTs, $R_{p,s} = V_{N-p,M-s} Q_{p,s}$.
3. The fast computation of the quaternion correlation of two images can be performed by three $N \times M$-point QDFTs. The number of quaternion multiplications is estimated as $3m_{2DDFT}$, or $6m_{2DDFT} + 16NM$ real multiplications.

- In the traditional (1,3)-model,
    1. The multiplicative property does not hold for the 2-D quaternion correlation, i.e., $R_{p,s} \neq V_{N-p,M-s} Q_{p,s}$.
    2. The computation of the correlation of two quaternion images requires at least $6m_{2DDFT} + 16NM$ real multiplications.
- The 2-D quaternion correlation functions in the (2,2)- and (1,3)-models are different functions.

The two-dimensional quaternion convolution in the (2,2) model has the same advantages over the (1,3) model as the correlation function described above.

|  | (1,3)-model | (2,2)-model |
|---|---|---|
| Multiplication | non-commutative | commutative |
| Correlation: | concept is not clear | unique concept |
|    1. Multiplicative property | not valid ($R_{p,s} \neq V_{N-p,M-s} Q_{p,s}$) | valid ($R_{p,s} = V_{N-p,M-s} Q_{p,s}$) |
|    2. Fast Algorithm | requires four right-side 2-D quaternion DFTs | requires six 2-D complex DFTs |
|    3. Number of real multiplications | $6m_{2DDFT} + 16NM$ | $\geq 6m_{2DDFT} + 16NM$ |

Table 1. Comparison of two models of quantum algebra.

## 5. CONCLUSION

An effective calculation method of the correlation function in the commutative (2,2)-model of quaternions is described for signals and color images. The correlation between quaternion signals and images is calculated by using the multiplication of their quaternion DFTs. The complexity of the correlation of color images is two times higher than in complex algebra. We plan to use (1) the developed tools in different practical applications and (2) random correlation concepts in the quaternion domain.

Conflicts of Interest: The author(s) declare that there are no conflicts of interest regarding the publication of this paper.